\documentclass{article}

\usepackage{graphicx}

\title{On the singular braid monoid of an orientable surface}
\author{Jer\'onimo D\'{\i}az--Cantos \\ Depto. de \'Algebra
\\ Universidad de Sevilla \and
Juan Gonz\'alez--Meneses\footnote{Supported by BFM 2001--3207 and FQM 218.}
\\ Depto. de Matem\'atica Aplicada I \\ Universidad de Sevilla \and
Jos\'e M. Tornero$^*$ \\ Depto. de \'Algebra \\ Universidad de Sevilla}
\date{August, 2002}

\newcommand{\df}{\noindent {\bf Definition.-- }}

\newcommand{\thr}{\noindent {\bf Theorem.-- }}

\newcommand{\dem}{{\bf Proof.-- }}
\newcommand{\obs}{\noindent {\bf Remark.-- }}
\newcommand{\cor}{\noindent {\bf Corollary.-- }}

\newcommand{\wh}{\widehat}
\newcommand{\ups}{\upsilon}

\begin{document}

\maketitle

\abstract{In this paper we show that the singular braid monoid of
an orientable surface can be embedded in a group. The proof is
purely topological, making no use of the monoid presentation.}

\vspace{.5cm}


\vspace{.5cm}


\section{Introduction}

The aim of this paper is proving, in an easy and fully topological
way, that the singular braid monoid of an orientable surface can
be embedded into a group.

The singular braid monoids appear naturally when studying braid
groups of surfaces. They were introduced in \cite{Baez} and \cite{B1}
and they arise in some situations connected with the Vassiliev
knot invariants (\cite{BN}).

A particular case of our theorem (that of the open disk) was proved by
Fenn, Keyman and Rourke in \cite{FKR}, where it was pointed out
that all the topological difficulties of the proof relied on the
so--called diamond lemma. Their proof, although purely
topological, was both quite involved and not adaptable (at least
in an easy way) for general orientable surfaces. On the other hand
our approach results in a simpler proof, even for the open disk case.

An algebraic generalization of this result can be found in
\cite{CIA}, which in fact may include the surface case if applied
together with the results of \cite{GM} concerning the monoid presentation.
However, we think that a topological and more down--to--earth proof of
this result may shed some light and contribute towards a better
understanding of the subject.

A final comment is in order : a different proof of this result, also based in
the monoid presentation, has been announced to us by Bellingeri in a
private communication.

\section{The singular braid monoid and the singular braid group}

\df Let $U$ be any orientable surface. If we fix beforehand $n$ distinct points $P_1,...,P_n \in U$,
a geometric braid is a set of $n$ differentiable disjoint paths $b_1,...,b_n$ inside the cylinder
$U \times [0,1]$, such that $b_i$ (also called the $i$--th string) goes, monotonically in
$t \in [0,1]$, from $\left( P_i,0 \right)$ to $\left( P_j,1 \right)$.

A singular geometric braid is just the same as a geometric braid, except for the fact that we allow
a finite number of transversal intersections between two strings. The intersection points will be
called singular points.

If we consider the braids modulo isotopies of $U \times [0,1]$ leaving fixed $U \times
\{ 0,1 \}$, we get the singular braid monoid (of $U$), noted $SB_n$,  with the operation
defined by concatenation.

We will note by $\alpha, \beta, \gamma,...$ the geometric braids and by
$\wh{\alpha},  \wh{\beta}, \wh{\gamma},...$ their corresponding classes in $SB_n$.

When $U$ is the open disk $D$, $SB_n$ is a very well--known monoid (\cite{B2}),
generated by the braids $ \wh{\sigma_i}, \wh{\sigma^{-1}_i}, \wh{\tau_i}$,
where $\sigma_i,\sigma_i^{-1},\tau_i$ are the geometric braids shown below.  We will also speak
of $ \wh{\sigma_i}, \wh{\sigma^{-1}_i}, \wh{\tau_i}$ in the general case, by means of
the natural embedding of the singular braid monoid of the disk into $SB_n$.

\begin{figure}[ht]
\centerline{\includegraphics{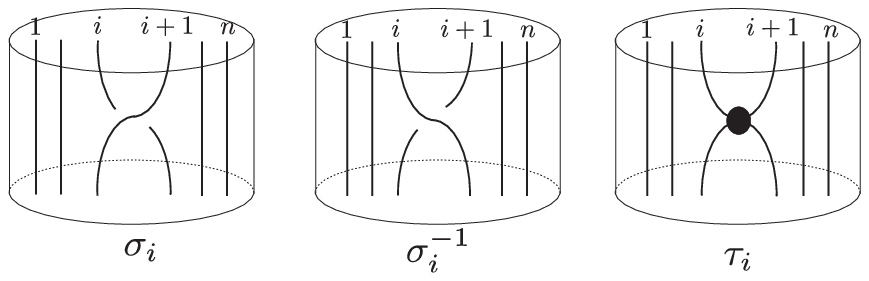}}
\end{figure}

In what follows, we will draw braids only in the cylinder $D \times [0,1]$ (for the sake
of simplicity), although our arguments will work for $U \times [0,1]$.

We will introduce now the analogous setup to that of Fenn, Keyman and Rourke in \cite{FKR}
for the disk case.

\df We will call $M$ the monoid where the elements are geometric singular braids (modulo
isotopies) in which every singular point is assigned a colour, black or white.

We will then note by $\alpha, \beta, \gamma,...$ the geometric braids (with coloured singular
points) and by $\wh{\alpha},  \wh{\beta}, \wh{\gamma},...$ their corresponding classes
in $M$.

\obs There is a surjection from $M$ to $SB_n$ which consists on forgetting the colours of the
singular points. On the other hand, $SB_n$  can be embedded in $M$ by assigning the
black colour to every singular point. The image of $\wh{\tau_i}$ under this injection
will be noted also by $\wh{\tau_i}$, while we will note by $\ups_i$ the braid obtained by
assigning the white colour to the singular point of $\tau_i$. We will call $\wh{\ups_i}$ the
opposite of $\wh{\tau_i}$.

\begin{figure}[ht]
\centerline{\includegraphics{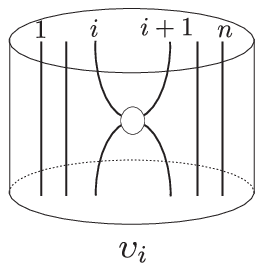}}
\end{figure}

Without explicit mention we will consider $SB_n$ as a submonoid of $M$, using the above
injection.

\df If we add to $M$ the relations $ \wh{\tau_i}\wh{\ups_i}
= \wh{\ups_i}\wh{\tau_i} =1$, we obtain a group, called the singular braid group
(of $U$), noted $SG_n$.

\begin{figure}[ht]
\centerline{\includegraphics{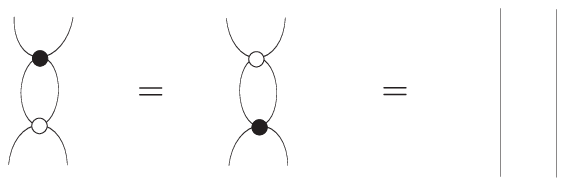}}
\centerline{The relations $ \wh{\tau_i}\wh{\ups_i}= \wh{\ups_i}\wh{\tau_i} = 1$ in $SG_n$.}
\end{figure}

From now on, if a singular braid has the form, say $\wh{\alpha} = \wh{\alpha_1 \tau_i \alpha_2}$, and
we do not care about who $\alpha_1$ and $\alpha_2$ are, we will draw it as follows:

\begin{figure}[ht]
\centerline{\includegraphics{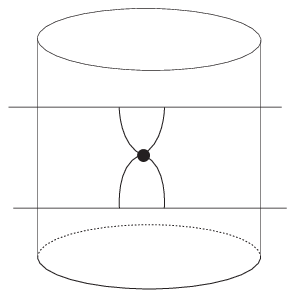}}
\end{figure}

\section{The embedding theorem}

\thr The natural map $SB_n \longrightarrow SG_n$ is one--to--one.

Before proving the result we will introduce a few notations. For $\alpha, \beta$ geometric
singular braids, we will note $\alpha \nearrow \beta$ if there exist $\alpha_1, \beta_1$
such that $\wh{\alpha} = \wh{\alpha_1}$, $\wh{\beta} = \wh{\beta_1}$ and we can obtain
$\beta_1$ from $\alpha_1$ by adding a pair of consecutive opposite singular points. In
the same way, we will note $\alpha \searrow \beta$ if there exist $\alpha_1, \beta_1$
such that $\wh{\alpha} = \wh{\alpha_1}$, $\wh{\beta} = \wh{\beta_1}$ and we can obtain
$\beta_1$ from $\alpha_1$ by erasing a pair of consecutive opposite singular points.

\obs Note that $\wh{\alpha}$ and $\wh{\beta}$ define the same element on $SG_n$ if and
only if there exist $\alpha = \alpha_0, \alpha_1,...,\alpha_k=\beta$ such that either $\alpha_i
\nearrow \alpha_{i+1}$ or $\alpha_i \searrow \alpha_{i+1}$ for all $i=0,...,k-1$.

The following result is the key part of the proof of the theorem, as noted in \cite{FKR}, where
it is proved for the case $U=D$. The proof will be given in the following section.

\noindent {\bf Diamond lemma.--} Let $\alpha, \beta, \gamma$ be geometric singular braids
such that $\alpha \nearrow \beta \searrow \gamma$. Then, either $\wh{\alpha} =
\wh{\gamma}$, or there exists $\eta$ such that $\alpha \searrow \eta \nearrow \gamma$.

\df We will say $\wh{\alpha} \in M$ is irreducible if there is no $\beta$
with $\alpha \searrow \beta$.

\cor If $\wh{\alpha}, \wh{\beta} \in M$ are irreducible and define the same element in
$SG_n$, then $\wh{\alpha} = \wh{\beta}$ in $M$.

The proof is immediate from the diamond lemma (see \cite{FKR}) and, besides, as the
elements of $SB_n$ are irreducible (they only have black singular points),
this proves that the map from $SB_n$ to $SG_n$ is an embedding. Hence the theorem
is proved, up to the diamond lemma.

\section{The diamond lemma}

Given a geometric singular braid $\beta$, with a singular point $p$ on it, we will denote
by $\beta (p^+)$ and $\beta(p^-)$ the braid obtained from $\beta$ by replacing $p$ by
a positive and negative crossing, respectively. That is, if $\beta = \beta_1 \tau_i \beta_2$,
(where $p$ is the singular point in $\tau_i$) we will have $\beta(p^+) =
\beta_1 \sigma_i \beta_2$ and $\beta(p^-) =  \beta_1 \sigma^{-1}_i \beta_2$. In a similar
way, if we have $p_1,...,p_m$ singular points on $\beta$, we will write
$\beta (p_1^{s_1},...,p_m^{s_m})$ the braid obtained by replacing each $p_i$ by a crossing
with sign $s_i \in \{+,-\}$.

If we have a singular braid $\beta$, with a singular point $p$ on it, we will write $\beta
(p^\bullet)$ and $\beta (p^\circ)$ the braid obtained from $\beta$ by replacing $p$ (no
matter which colour it has) by a black point and a white point, respectively.

\noindent {\bf Lemma 1.--} Let $\beta$ and $\beta'$ be geometric singular braids, $p$
a singular point on $\beta$. Assume $\wh{\beta} = \wh{\beta'}$ and let us call also $p$
the point corresponding to $p$ in $\beta'$. Then:
\begin{enumerate}
\item $\wh{\beta(p^+)} = \wh{\beta'(p^+)}$.
\item $\wh{\beta(p^-)} = \wh{\beta'(p^-)}$.
\end{enumerate}

\dem As both cases are analogous, we will do the first case only. As $\wh{\beta} =
\wh{\beta'}$ , there exists an isotopy of the cylinder, $H_t$, with $H_0 (\beta) = \beta$,
$H_1 (\beta) = \beta'$.
Take a sphere $S$ centered at $p$ with radius small enough such that the only strings
intersecting $S$ are those forming the singular point. We can assume that $\beta(p^+)$ and
$\beta$ coincide outside $S$.

\begin{figure}[ht]
\centerline{\includegraphics{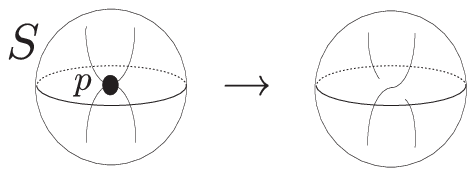}}
\end{figure}

We can suppose that $H_1(S)$ is also a sphere centered at
$p$ on $\beta'$. Hence, if we apply $H_1$ to $\beta(p^+)$, we obtain the braid $\beta'$,
except for the fact that the point
$p$ has been replaced by a positive crossing (since $U$ is orientable). That is, we get
$\beta'(p^+)$. So $H_1(\beta(p^+))= \beta'(p^+)$.

\noindent {\bf Lemma 2.--} Let $\wh{\beta} \in M$, having two consecutive singular
 points, $p$ and $q$. Then
$
\wh{\beta(p^s,q^c)} = \wh{\beta(p^c,q^s)}, \mbox{ with } s \in \{+,-\}, \; c \in \{
\bullet, \circ \}.
$

\dem This result is a straightforward consequence of the well known relation
$\sigma_i \tau_i = \tau_i \sigma_i$ in the
singular braid monoid of the disk.

\begin{figure}[ht]
\centerline{\includegraphics{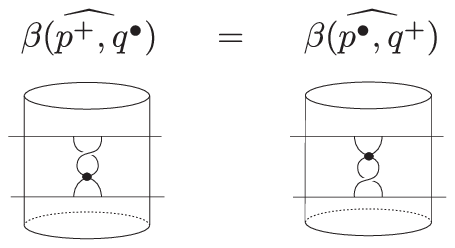}}
\end{figure}

We can now proceed to prove the diamond lemma. As $\alpha \nearrow \beta$
we can consider, with no loss of generality, that $\alpha= \beta(p^+,q^-) $ for some pair
$p,q$ of opposite consecutive points in $\beta$. On the other side, as $\beta \searrow
\gamma$, let us call $\beta'$ the braid verifying $\wh{\beta} = \wh{\beta'}$ and, as above,
$\gamma = \beta'(r^+,s^-)$ for a pair $r,s$ of opposite consecutive points in $\beta'$.

Let us call $H$ the cylinder isotopy taking $\beta$ into $\beta'$. Then, using our
previous notation, we can write
$
\alpha = \beta(p^+,q^-) \nearrow \beta \stackrel{H}{\longrightarrow} \beta' \searrow
\beta'(r^+,s^-) = \gamma.
$

The technique used for the proof is heavily related to the proof of lemma 1. In fact, we will
make an extensive use of the fact that $H$ brings together $r$ and $s$, while possibly
moving apart
$p$ and $q$, and $H^{-1}$ acts the other way around. This fact, together with the using of
the small spheres as in lemma 1, will give us the result.

We will assume that $p$ is above $q$ and $r$ is above $s$ from now on. Then we need
to distinguish three cases:

\begin{enumerate}
\item[(a)] $p=r$, $q=s$.
\item[(b)] $p, q, r, s$ are all distinct.
\item[(c)] $p=s$ or $q=r$.
\end{enumerate}

\noindent {\bf Case (a)} This is the easiest situation. As $\wh{\beta} = \wh{\beta'}$, from
lemma 1 we  have $\wh{\beta(p^+)} =  \wh{\beta'(p^+)}$. Hence
$\wh{\alpha} = \wh{\beta(p^+,q^-)} =  \wh{\beta'(p^+,q^-)} = \wh{\gamma}$.

\noindent {\bf Case (b)} In this case $\beta$ has two pairs of opposite singular points, $(p,q)$
and $(r,s)$, of which $(p,q)$ are consecutive. Hence $\beta$ can be assumed to have the
form $\beta_1 \tau_i \ups_i \beta_2$, and the analogous thing happens with $\beta'$ and
$(r,s)$.

\begin{figure}[ht]
\centerline{\includegraphics{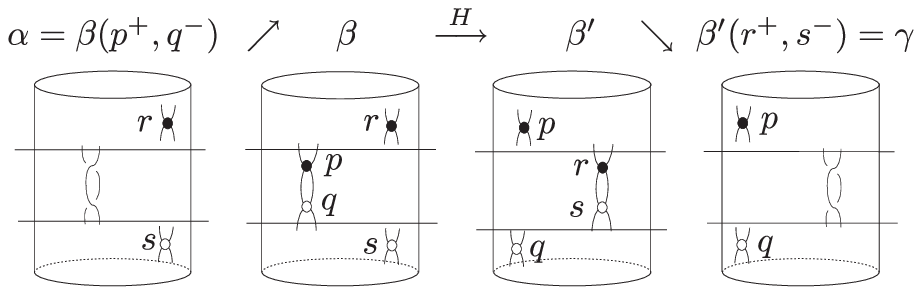}}
\end{figure}

Then we can write
$$
\alpha = \beta(p^+,q^-) \stackrel{H}{\longrightarrow} \beta'(p^+,q^-) \searrow
\beta'(p^+,q^-,r^+,s^-) \stackrel{H^{-1}}{\longrightarrow}
$$
$$
\stackrel{H^{-1}}{\longrightarrow} \beta(p^+,q^-,r^+,s^-) \nearrow
\beta(r^+,s^-) \stackrel{H}{\longrightarrow} \beta'(r^+,s^-) = \gamma.
$$

\begin{figure}[ht]
\centerline{\includegraphics{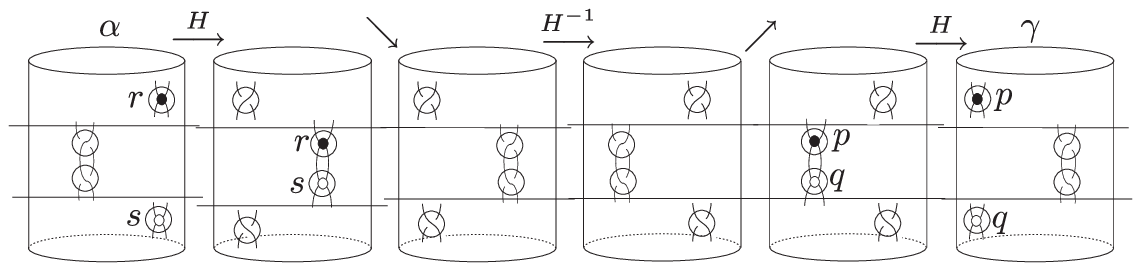}}
\end{figure}

This proves the existence of $\eta = \beta'(p^+,q^-,r^+,s^-) $ with $\alpha \searrow \eta
\nearrow \gamma$.

\noindent {\bf Case (c)} We will do the subcase $q=r$, the other one being symmetric. We
have

\begin{figure}[ht]
\centerline{\includegraphics{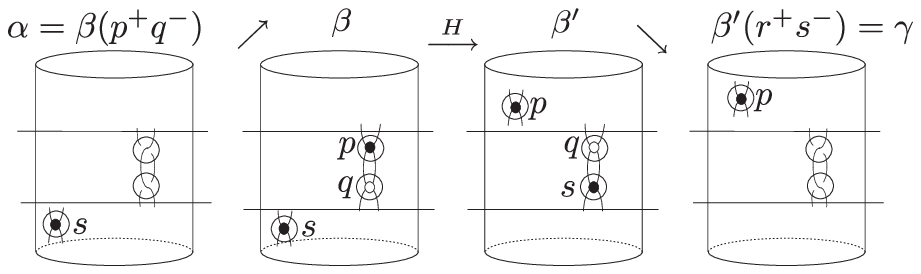}}
\end{figure}

Then, from lemma 2,
$$
\alpha = \beta(p^+,q^-) \stackrel{H}{\longrightarrow} \beta'(p^+,q^-) =
\beta'(p^+,q^-,s^\bullet) \stackrel{\sim}{\longrightarrow} \beta'(p^+,q^\bullet,s^-)
\stackrel{H^{-1}}{\longrightarrow}
$$
$$
\stackrel{H^{-1}}{\longrightarrow} \beta(p^+,q^\bullet,s^-) \stackrel{\sim}{\longrightarrow}
\beta(p^\bullet,q^+,s^-) \stackrel{H}{\longrightarrow} \beta'(p^\bullet,q^+,s^-) = \gamma,
$$
where $\stackrel{\sim}{\longrightarrow}$ stands for the (non--specified) isotopies which
come from applying the lemma 2. This proves $\wh{\alpha} = \wh{\gamma}$ and concludes
the proof of the diamond lemma.

\begin{figure}[ht]
\centerline{\includegraphics{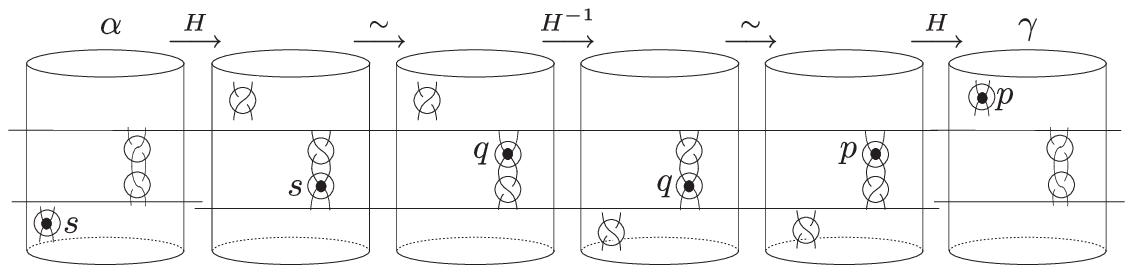}}
\end{figure}


\begin{thebibliography}{99}

\bibitem[\bf 1]{Baez}
{\sc J. Baez,} `Link invariants and perturbation theory',
{\em Lett. Math. Phys.}  2 (1992) 43--51.

\bibitem[\bf 2]{BN}
{\sc D. Bar--Natan,} `On the Vassiliev knot invariants',
{\em Topology} (2) 34 (1995) 423--472.

\bibitem[\bf 3]{CIA}
{\sc G. Basset,} `Quasi--commuting extensions of groups',
{\em Comm. Algebra} (11) 28 (2000) 5443--5454.

\bibitem[\bf 4]{B1}
{\sc J. Birman,} {\em Braids, links and mapping class groups}.
(Princeton University Press, Princeton, 1974).

\bibitem[\bf 5]{B2}
{\sc J. Birman,} `New points of view in knot theory', {\em B.
Am. Math. Soc.} (2) 28 (1993) 253--287.

\bibitem[\bf 6]{FKR}
{\sc R. Fenn, E. Keyman, C. Rourke,} `The singular braid
monoid embeds in a group', {\em J. Knot Theor. Ramif.}
(7) 7 (1998) 881--892.

\bibitem[\bf 7]{GM}
{\sc J. Gonz\'alez--Meneses,} `Presentations for the monoids
of singular braids on closed surfaces', {\em Comm. Algebra}
(6) 30 (2002) 2829--2836.

\end{thebibliography}
\end{document}